\documentclass[12pt]{article}

\usepackage[margin=1in]{geometry}  
\usepackage{graphicx}              
\usepackage{amsmath}               
\usepackage{amsfonts}              
\usepackage{amsthm}                
\usepackage{amssymb}
\usepackage{amscd}
\usepackage{amsmath,amscd}
\usepackage{young}

\begin{document}

\nocite{*}

\title{Poisson cohomology of two Fano threefolds.}

\author{Evgeny Mayanskiy}

\maketitle

\begin{abstract}
  We study the variety of Poisson structures and compute Poisson cohomology for two families of Fano threefolds - smooth cubic threefolds and the del Pezzo quintic threefold. Along the way we reobtain by a different method earlier results of Loray, Pereira and Touzet in the special case we are considering. 
\end{abstract}

\setcounter{section}{-1}
\section{Introduction.}

A {\it Poisson structure} on a smooth algebraic variety $X$ is given by a bivector field $\omega\in H^0(X,{\wedge}^2 T_X)$, which satisfies the condition
$$
[\omega, \omega]=0, \;\;\;\;\;\;\;\;\;\; (\star )
$$
where $[\; , \; ]$ denotes the Schouten bracket.\\

It is natural to consider the variety of Poisson structures $\mathcal P$ on $X$, which is the subvariety in ${\mathbb P}(H^0 (X,{\wedge}^2T_X))$ given by the equations $(\star )$.\\

Given any point $\omega \in \mathcal P$, one defines Poisson cohomology $H_{Poisson}^\star (X,\omega)$ of $X$ with respect to $\omega$ \cite{Lichnerowicz}.\\

Explicit computation of Poisson cohomology for various algebraic varieties is of interest (see \cite{HongXu}). In \cite{HongXu} Wei Hong and Ping Xu computed Poisson cohomology of del Pezzo surfaces. The condition $(\star )$ is automatically satisfied in that case for dimension reasons. This is not the case for Fano threefolds.\\

Poisson Fano threefolds with Picard number $1$ were classified in \cite{Loray1} and \cite{Loray2}. They are exactly the following (see \cite{Loray1}, \cite{Loray2}):
\begin{itemize}
\item the projective space ${\mathbb P}^3$,
\item the quadric $Q\subset {\mathbb P}^4$,
\item a sextic hypersurface $X\subset {\mathbb P} (1,1,1,2,3)$,
\item a quartic hypersurface $X\subset {\mathbb P} (1,1,1,1,2)$,
\item a cubic hypersurface $X\subset {\mathbb P}^4$,
\item a complete intersection of two quadrics in ${\mathbb P}^5$,
\item the del Pezzo quintic threefold $X\subset {\mathbb P}^6$,
\item the Mukai-Umemura threefold.
\end{itemize}

In each case Loray, Pereira and Touzet described the dimensions of the irreducible components of the variety of Poisson structures, described their smoothness and some other properties (see \cite{Loray1}, \cite{Loray2}).\\

In this note we consider (from a different viewpoint) two families of Poisson Fano threefolds from the above list: cubic threefolds and the del Pezzo quintic threefold. For each of them we describe explicitly the variety of Poisson structures (thus reobtaining the results of \cite{Loray1}, \cite{Loray2} in a special case) and compute Poisson cohomology.\\

Our main results are the following. For the explicit computations of matrices $A_{\omega}$, $B_{\omega}$ and $C_{\omega}$ we refer the reader to sections 2.3 and 3.2.\\

{\bf Theorem 1 (Loray, Pereira, Touzet, \cite{Loray1}, \cite{Loray2}).} {\it Let $X$ be a smooth cubic threefold. Then the variety of Poisson structures $\mathcal P\subset {\mathbb P}(H^0(X,{\wedge}^2T_X))$ on $X$ is the Grassmannian $G(2,5)$ of lines in ${\mathbb P}^4$ and $\mathcal P\subset {\mathbb P}(H^0(X,{\wedge}^2T_X))$ is the Pl{\" u}cker embedding.}\\

{\bf Theorem 2.} {\it Let $X$ be a smooth cubic threefold and $\omega\in \mathcal P\subset {\mathbb P}(H^0(X,{\wedge}^2T_X))$ a Poisson structure on $X$. Then the dimensions of the Poisson cohomology groups are as follows:
\begin{itemize}
\item $dim(H^0_{Poisson}(X,\omega))=1$,
\item $dim(H^1_{Poisson}(X,\omega))=0$,
\item $dim(H^2_{Poisson}(X,\omega))=20-rk(C_{\omega})$,
\item $dim(H^3_{Poisson}(X,\omega))=15-rk(C_{\omega})$.
\end{itemize}}

{\bf Theorem 3 (Loray, Pereira, Touzet, \cite{Loray1}, \cite{Loray2}).} {\it Let $X$ be the (smooth) del Pezzo quintic threefold. Then the variety of Poisson structures $\mathcal P\subset {\mathbb P}(H^0(X,{\wedge}^2T_X))$ on $X$ is the disjoint union of the Grassmannian $G(2,7)$ of lines in ${\mathbb P}^6$ (embedded into ${\mathbb P}(H^0(X,{\wedge}^2T_X))$ by the Pl{\" u}cker embedding) and a smooth conic. The plane spanned by the conic does not intersect the Grassmannian.}\\

{\bf Theorem 4.} {\it Let $X$ be the (smooth) del Pezzo quintic threefold and $\omega\in \mathcal P\subset {\mathbb P}(H^0(X,{\wedge}^2T_X))$ a Poisson structure on $X$. Then the dimensions of the Poisson cohomology groups are as follows:
\begin{itemize}
\item $dim(H^0_{Poisson}(X,\omega))=1$,
\item $dim(H^1_{Poisson}(X,\omega))=3-rk(A_{\omega})$,
\item $dim(H^2_{Poisson}(X,\omega))=21-rk(A_{\omega})-rk(B_{\omega})$,
\item $dim(H^3_{Poisson}(X,\omega))=23-rk(B_{\omega})$.
\end{itemize}}

The identification of the irreducible components of the varieties of Poisson structures with Grassmannians was proved in greater generality in \cite{Loray1}, \cite{Loray2}. In the cases we are considering we reobtain these results with a different method. Thus Theorem 1 and Theorem 3 is entirely a result of Loray, Pereira and Touzet. Note that in the case of the del Pezzo quintic threefold the curve component is shown to be smooth and rational in \cite{Loray1}, \cite{Loray2}. The fact that it is a conic was not mentioned explicitely there, but Loray, Pereira and Touzet knew it thanks to the joint work of Jorge Pereira and Carlo Perrone \cite{Pereira}. We thank Jorge Pereira for the last comment as well as for the following Remark.\\

{\bf Remark.} Still, note that the arguments of \cite{Loray1}, \cite{Loray2} do not shed light on the scheme structure of the space of Poisson structures.\\

The key ingredient of our computation of Poisson cohomology for cubics and the del Pezzo quintic threefold is the same as in \cite{HongXu} - the spectral sequence of Laurent-Gengoux, Sti{\'e}non and Xu \cite{Stienon}.\\

We work over $k=\mathbb C$ throughout.\\

\section{Schouten bracket.}

The Schouten bracket on a smooth algebraic variety $X$ gives a way to extend the structure of a Lie algebra on vector fields on $X$ to the structure of a graded Lie algebra on all multivector fields on $X$.\\

A convenient local definition of the Schouten bracket is given in a preprint by Bondal \cite{Bondal}:\\

{\bf Definition (Bondal, \cite{Bondal}). }  Let $A$, $B$ be multivector fields on $X$ of degrees $n$ and $m$ respectively. Then their Schouten bracket $[A,B]$ is a multivector field of degree $n+m-1$ such that for any $(n+m-1)-$form $\omega$ we have
$$
[A,B](\omega)=(-1)^{m(n-1)}A(dB(\omega))+(-1)^{n}B(dA(\omega))-AB(d\omega ).
$$

When $dim(X)=3$, one can use the isomorphism of vector bundles ${\wedge}^2T_X\cong {\Omega}^1_X(-K_X)$ and reformulate this definition as follows.\\

{\bf Proposition 1.} {\it Assume that $X$ is a smooth projective algebraic threefold. Let ${\omega}_A, {\omega}_B \in H^0(X,{\Omega}^1_X(-K_X))$ be the forms corresponding to bivector fields $A,B\in H^0(X,{\wedge}^2T_X)$ under the isomorphism ${\wedge}^2T_X\cong {\Omega}^1_X(-K_X)$. Then the Schouten bracket $[{\omega}_A, {\omega}_B]$ is equal to
$$
[{\omega}_A, {\omega}_B]=\frac{1}{Vol}\cdot \left( {\omega}_A \wedge d{\omega}_B + d{\omega}_A \wedge {\omega}_B \right).
$$
}\\

Note that in this case the Schouten bracket 
$$
[\; , \; ]\colon H^0(X,{\wedge}^2T_X) \times H^0(X,{\wedge}^2T_X)\rightarrow H^0(X,{\wedge}^3T_X)\cong H^0(X,-K_X)
$$
has its image in the space of global sections of the anticanonical line bundle $-K_X$, and the division by $Vol$ in the formula above signifies the isomorphism ${\Omega}^3_X\cong K_X$. In other words, if $dim(X)=3$, then in terms of forms the Schouten bracket
$$
[\; , \; ]\colon H^0(X,{\Omega}^1_X(-K_X)) \times H^0(X,{\Omega}^1_X(-K_X))\rightarrow H^0(X,-K_X)
$$
is the composition of the exterior derivatives and the exterior product:
\begin{multline*}
H^0(X,{\Omega}^1_X(-K_X)) \times H^0(X,{\Omega}^1_X(-K_X))\xrightarrow{1\times d+d\times 1} [ H^0(X,{\Omega}^1_X(-K_X)) \times H^0(X,{\Omega}^2_X(-K_X)) ] \oplus \\
[ H^0(X,{\Omega}^2_X(-K_X)) \times H^0(X,{\Omega}^1_X(-K_X)) ] \xrightarrow{\wedge} H^0(X,{\Omega}^3_X(-2K_X))\cong H^0(X,-K_X).
\end{multline*}

{\it Proof:} This follows from Bondal's definition by a local computation.\\

Let $x_1,x_2,x_3$ be local coordinates on $X$. As an element of $H^0(X,{\wedge}^2T_X)$, ${\omega}_A$ maps $dx_i\wedge dx_j$ to $\frac{1}{Vol}\cdot {\omega}_A \wedge dx_i\wedge dx_j$, where $Vol=dx_1\wedge dx_2 \wedge dx_3$.\\

Hence 
\begin{multline*}
[A,B](dx_1\wedge dx_2 \wedge dx_3)=A\left( d\frac{1}{Vol}({\omega}_B \wedge dx_1\wedge dx_2)\wedge dx_3 - d \frac{1}{Vol} ({\omega}_B \wedge dx_1\wedge dx_3)\wedge dx_2 +\right. \\
\left. +d\frac{1}{Vol} ({\omega}_B \wedge dx_2\wedge dx_3)\wedge dx_1   \right) +(A\leftrightarrow B)=\\
=\frac{1}{Vol}\cdot {\omega}_A \wedge [dc_1 \wedge dx_1 + dc_2 \wedge dx_2 + dc_3 \wedge dx_3]+(A\leftrightarrow B),
\end{multline*}
where $c_i$ are such that ${\omega}_B=\sum_{i=1}^{3} c_i\cdot dx_i$.\\

This expression is exactly equal to
$$
\frac{1}{Vol} {\omega}_A \wedge d{\omega}_B + (A\leftrightarrow B).
$$
{\it QED}\\

If $X$ is a Fano threefold of index $2$ (the only case we are considering), then $K_X\cong {\mathcal O}_X(-2)$. In this case, the Schouten bracket is a bilinear map
$$
[\; , \; ]\colon H^0(X,{\Omega}^1_X(2))\times H^0(X,{\Omega}^1_X(2))\rightarrow H^0(X,{\mathcal O}_X(2)).
$$

From now on $X$ will always denote a smooth Fano threefold of index $2$. We will always identify $H^0(X,{\wedge}^2T_X)$ and $H^0(X,{\Omega}^1_X(2))$ (as well as $H^0(X,-K_X)$ and $H^0(X,{\mathcal O}(2))$) by the isomorphism ${\wedge}^2T_X\cong {\Omega}^1_X(-K_X)$.\\

\section{Cubic threefolds.}

Let $X\subset {\mathbb P}^4$ be a smooth cubic threefold. Let us denote by $Z_0,Z_1,Z_2,Z_3,Z_4$ the homogeneous coordinates in ${\mathbb P}^4$ and by
$$
F=F(Z_0,Z_1,Z_2,Z_3,Z_4)=0
$$
the equation of $X$.\\

\subsection{Cohomology computations.}

In order to describe the variety of Poisson structures on $X$ and compute Poisson cohomology of $X$, we need to find handy descriptions for the spaces of multivector fields $H^0(X,{\wedge}^jT_X)$ on $X$ as well as to compute (or check the vanishing) of the higher cohomology groups $H^i(X,{\wedge}^jT_X), i\geq 1$. This is the subject of the present subsection. The methods we are using are standard and the results are well-known.\\

{\bf Lemma 1.} {\it Let $X\subset {\mathbb P}^4$ be a smooth cubic threefold. Then 
\begin{itemize}
\item[(a)] $H^0(X,{\mathcal O}_X)=k$, $H^i(X,{\mathcal O}_X)=0$ for any $i\geq 1$,
\item[(b)] $H^i(X,T_X)=0$ for any $i\neq 1$,
\item[(c)] $H^1(X,T_X)$ is a $10$-dimensional vector space, which can be identified with the cokernel of the homomorphism
$$
H^0({\mathbb P}^4,\mathcal O)\oplus H^0({\mathbb P}^4,\mathcal O (1))^{\oplus 5} \rightarrow H^0({\mathbb P}^4,{\mathcal O}(3)), (a,(b_0,b_1,b_2,b_3,b_4))\mapsto aF+ \sum_{i=0}^{4} b_i \frac{\partial F}{\partial Z_i},
$$
\item[(d)] $H^i(X,{\wedge}^2T_X)=0$ for any $i\geq  1$,
\item[(e)] $H^0(X,{\wedge}^2T_X)$ is a $10$-dimensional vector space, which can be identified with the kernel of the (surjective) homomorphism
$$
H^0({\mathbb P}^4,\mathcal O (1))^{\oplus 5} \rightarrow H^0({\mathbb P}^4,{\mathcal O}(2)), (a_0,a_1,a_2,a_3,a_4)\mapsto \sum_{i=0}^{4} a_i Z_i,
$$
\item[(f)] $H^i(X,{\wedge}^3T_X)=0$ for any $i\geq  1$,
\item[(g)] $H^0(X,{\wedge}^3T_X)\cong H^0({\mathbb P}^4,{\mathcal O}(2))$ has dimension $15$.
\end{itemize}}

{\it Proof:}\\ 
(a) The short exact sequence of sheaves on ${\mathbb P}^4$
\begin{equation}
\label{SES1}
0\rightarrow {\mathcal O}(-3)\rightarrow \mathcal O \rightarrow {\mathcal O}_X \rightarrow 0 
\end{equation}

gives the long exact sequence of cohomology groups 
$$
0\rightarrow H^0({\mathbb P}^4,{\mathcal O}(-3))\rightarrow  H^0({\mathbb P}^4,O) \rightarrow H^0(X,{\mathcal O}_X) \rightarrow H^1({\mathbb P}^4, {\mathcal O}(-3))\rightarrow \ldots
$$

Since $H^i({\mathbb P}^4, {\mathcal O}(-3))=0$ for any $i\in \mathbb Z$, we obtain the isomorphisms
$$
H^i({\mathbb P}^4,O) \cong H^i(X,{\mathcal O}_X),\;\;\; i\in \mathbb Z.
$$
Since $H^i({\mathbb P}^4, {\mathcal O})=0$ for any $i\geq 1$, we conclude.\\

(b)(c) Multiplying the exact sequence~(\ref{SES1}) by ${\mathcal O}(d)$, $d\in\mathbb Z$, we get
\begin{equation}
\label{SES2}
0\rightarrow {\mathcal O}(d-3)\rightarrow \mathcal O (d) \rightarrow {\mathcal O}_X (d) \rightarrow 0. 
\end{equation}

Since $H^i({\mathbb P}^4, {\mathcal O} (d))=0$ for $0< i < 4$, $d\in \mathbb Z$ and for $i=4, d\geq -4$, the long exact sequence of cohomology groups implies that 
\begin{gather*}
H^i(X, {\mathcal O}_X (d))=0 \;\mbox{for }\; 0<i<3, d\in \mathbb Z,\\
H^3(X, {\mathcal O}_X (d))=0 \;\mbox{for }\; d\geq -1,\\
H^4(X, {\mathcal O}_X (d))=0 \;\mbox{for }\; d\geq -4.\\
\end{gather*}

Moreover, we get the short exact sequences 
\begin{equation}
\label{SES3}
0\rightarrow H^0({\mathbb P}^4,{\mathcal O}(d-3)) \rightarrow H^0({\mathbb P}^4,{\mathcal O}(d)) \rightarrow H^0(X,{\mathcal O}_X(d)) \rightarrow 0, 
\end{equation}
which compute the spaces of global sections of line bundles ${\mathcal O}_X(d)$ on $X$ in terms of those on ${\mathbb P}^4$.\\

In particular, for $d\leq 2$ we get isomorphisms
$$
H^0({\mathbb P}^4,{\mathcal O}(d))\cong H^0(X,{\mathcal O}_X(d)).
$$

Let $j\colon X\hookrightarrow {\mathbb P}^4$ denote the embedding. Consider the normal bundle sequence on $X$
\begin{equation}
\label{SES5}
0\rightarrow T_X\rightarrow j^{\star}T_{{\mathbb P}^4} \rightarrow {\mathcal O}_X (3) \rightarrow 0 
\end{equation}
and the pullback to $X$ of the Euler exact sequence on ${\mathbb P}^4$
\begin{equation}
\label{SES6}
0\rightarrow {\mathcal O}_X\rightarrow {\mathcal O}_X (1)^{\oplus 5} \rightarrow j^{\star}T_{{\mathbb P}^4} \rightarrow 0. 
\end{equation}

Since $H^i(X,{\mathcal O}_X(d))=0$ for any $i\geq 1$ and any $d\geq 0$, we conclude that 
$$
H^i(X,j^{\star}T_{{\mathbb P}^4})=0 \;\;\; \mbox{for any}\; i\geq 1.
$$

It is shown in \cite{KodairaSpencer} that $H^0(X,T_X)=0$. Then the long exact sequence of cohomology groups of~(\ref{SES5}) gives that 
$$
H^i(X,T_X)=0, i\geq 2
$$
and also the short exact sequence
\begin{equation}
\label{SES7}
0\rightarrow H^0(X,j^{\star}T_{{\mathbb P}^4}) \rightarrow H^0(X,{\mathcal O}(3)) \rightarrow H^1(X,T_X) \rightarrow 0. 
\end{equation}

The space of global sections of $j^{\star}T_{{\mathbb P}^4}$ on $X$ can be computed from~(\ref{SES6}) as follows:
$$
0\rightarrow H^0(X,{\mathcal O}_X) \rightarrow H^0(X,{\mathcal O}_X(1))^{\oplus 5} \rightarrow H^0(X,j^{\star}T_{{\mathbb P}^4}) \rightarrow 0. 
$$

Since $H^0(X,{\mathcal O}_X) = H^0({\mathbb P}^4,{\mathcal O})$ and $H^0(X,{\mathcal O}_X (1)) \cong  H^0({\mathbb P}^4,{\mathcal O} (1))$, using Euler's exact sequence again we conclude that
$$
H^0(X,j^{\star}T_{{\mathbb P}^4}) \cong  H^0({\mathbb P}^4,T_{{\mathbb P}^4})
$$

is the cokernel of the following homomorphism of vector spaces
$$
H^0({\mathbb P}^4,{\mathcal O}) \rightarrow H^0({\mathbb P}^4,{\mathcal O}(1))^{\oplus 5}, 1\mapsto (Z_0,Z_1,Z_2,Z_3,Z_4).
$$

The long exact sequence of~(\ref{SES3}) computes $H^0(X,{\mathcal O}_X(3))$ as follows:
$$
0\rightarrow H^0({\mathbb P}^4,{\mathcal O}) \rightarrow H^0({\mathbb P}^4,{\mathcal O}(3)) \rightarrow H^0(X,{\mathcal O}_X(3)) \rightarrow 0. 
$$

The homomorphism $H^0(X,j^{\star}T_{{\mathbb P}^4}) \rightarrow H^0(X,{\mathcal O}_X(3))$ lifts to
$$
H^0({\mathbb P}^4,{\mathcal O}(1))^{\oplus 5} \rightarrow H^0({\mathbb P}^4,{\mathcal O}(3)), (a_0,a_1,a_2,a_3,a_4)\mapsto \sum_{i=0}^{4} a_i\cdot  \frac{\partial F}{\partial Z_i}.
$$

Indeed, the element of $H^0(X,j^{\star}T_{{\mathbb P}^4})$ corresponding to $(a_0,a_1,a_2,a_3,a_4)\in H^0({\mathbb P}^4, {\mathcal O}(1))^{\oplus 5}$ is the vector field $\sum_{i=0}^{4} a_i\cdot \frac{\partial }{\partial Z_i} {\mid}_X$. Its image in $H^0(X,{\mathcal O}_X(3))$ is the corresponding section of the normal bundle $N_{X / {\mathbb P}^4}\cong {\mathcal O}_X(3)$ to $X$ in ${\mathbb P}^4$.\\

This proves part (c).\\

(d)(e) Take the dual of the normal bundle sequence~(\ref{SES5}):
\begin{equation}
\label{SES8}
0\rightarrow {\mathcal O}(-3) \rightarrow j^{\star}\Omega_{{\mathbb P}^4}^1 \rightarrow \Omega_{X}^1 \rightarrow 0. 
\end{equation}

The long exact sequence of cohomology groups of~(\ref{SES2}) implies that $H^0(X,{\mathcal O}_X(d))=0$ for $d\leq -1$. Note that $H^i(X,{\mathcal O}_X(-1))=0$ for any $i\in \mathbb Z$. Hence~(\ref{SES8}) gives (after multiplying by $\mathcal O (2)$ and taking cohomology)
$$
H^i(X,j^{\star}\Omega_{{\mathbb P}^4}^1(2))\cong H^i(X,\Omega_{X}^1(2))\; \mbox{for any }\; i\in \mathbb Z.
$$

Taking the dual of the short exact sequence~(\ref{SES6}) and multiplying the result by $\mathcal O (2)$, we find
\begin{equation}
\label{SES9}
0 \rightarrow j^{\star}\Omega_{{\mathbb P}^4}^1(2) \rightarrow {\mathcal O}_{X}(1)^{\oplus 5}\rightarrow {\mathcal O}_X(2) \rightarrow 0. 
\end{equation}

Since $H^i(X,{\mathcal O}_X(1))=H^i(X,{\mathcal O}_X(2))=0,\;\; i\geq 1$, this implies immediately that 
$$
H^i(X,j^{\star}\Omega_{{\mathbb P}^4}^1(2))=0,\;\;\; i\geq 2.
$$

Moreover, by taking cohomology in~(\ref{SES9}) we obtain an exact sequence
$$
0\rightarrow H^0(X,j^{\star}\Omega_{{\mathbb P}^4}^1(2)) \rightarrow H^0(X,{\mathcal O}_X(1))^{\oplus 5}\rightarrow  H^0(X,{\mathcal O}_X(2))\rightarrow H^1(X,j^{\star}\Omega_{{\mathbb P}^4}^1(2) )\rightarrow 0.
$$

The homomorphism in the middle can be identified with
$$
H^0({\mathbb P}^4,\mathcal O (1))^{\oplus 5} \rightarrow H^0({\mathbb P}^4,{\mathcal O}(2)), (a_0,a_1,a_2,a_3,a_4)\mapsto \sum_{i=0}^{4} a_i Z_i.
$$

It is clearly surjective. Hence $H^1(X,j^{\star}\Omega_{{\mathbb P}^4}^1(2))=0$.\\

This gives the vanishing $H^i(X,{\wedge}^2T_X)\cong H^i(X,{\Omega}^1_X(2))\cong H^i(X,j^{\star}\Omega_{{\mathbb P}^4}^1(2))=0$ for any $i\geq 1$ as well as the description of $H^0(X,{\wedge}^2T_X)$ in part (e).\\

(f)(g) Note that ${\wedge}^3T_X\cong ({\Omega}^3_X)^{\vee}\cong {\mathcal O}(-K_X)\cong {\mathcal O}_X(2)$ since $X$ has index $2$.\\

Hence $H^i(X, {\wedge}^3T_X)\cong H^i(X, {\mathcal O}_X(2))\cong H^i({\mathbb P}^4, {\mathcal O}(2))$ for any $i\in \mathbb Z$. {\it QED}\\

{\bf Corollary 1.} {\it Let $X$ be a smooth cubic threefold. Then $H^0(X,{\wedge}^2T_X)$ can be identified (as a vector space) with $\mathfrak{so}(5)$, which is the Lie algebra of skew-symmetric $5\times 5$ matrices.\\

As a basis of $H^0(X,{\wedge}^2T_X) \cong H^0(X,{\Omega}_X^1(2))$ one can take $1$-forms 
$$
{\epsilon}_{ij}=Z_jdZ_i-Z_idZ_j, \; i<j.
$$
}\\

{\it Proof:} In Lemma 1(e) we described $H^0(X,{\wedge}^2T_X) \cong H^0(X,{\Omega}_X^1(2))$ as the kernel of the homomorphism
$$
H^0({\mathbb P}^4,\mathcal O (1))^{\oplus 5} \rightarrow H^0({\mathbb P}^4,{\mathcal O}(2)),\;\;\; (a_0,a_1,a_2,a_3,a_4)\mapsto \sum_{i=0}^{4} a_i\cdot Z_i.
$$

As a basis of this kernel one can take sequences of the form
$$
{\epsilon}_{ij}=(0,\cdots , Z_j,\cdots , -Z_i ,\cdots ,0), \; i<j,
$$
where we have $Z_j$ in the $i$-th position, $-Z_i$ in the $j$-th position and $0$ in the remaining positions.\\

Then an arbitrary element of the kernel will have the form $\sum_{0\leq i < j \leq 4} {\beta}_{ij}\cdot {\epsilon}_{ij}$, where $({\beta}_{ij})\in \mathfrak{so}(5)$.\\

A sequence $(a_0,a_1,a_2,a_3,a_4)\in H^0({\mathbb P}^4,\mathcal O (1))^{\oplus 5}$ such that $\sum_{i=0}^{4} a_i\cdot Z_i=0$ corresponds to the form $\sum_{i=0}^{4} a_i \cdot dZ_i$ on $X$. Hence ${\epsilon}_{ij}$ are exactly the forms $Z_jdZ_i-Z_idZ_j$. {\it QED}\\

More invariantly, this leads to an identification $H^0(X,{\wedge}^2T_X)\cong {\wedge}^2V$, where $V=H^0({\mathbb P}^4, {\mathcal O}(1))$ is the $5$-dimensional vector space such that $X\subset {\mathbb P}(V)$ is the vanishing locus of $F\in Sym^3V$.\\

We will see later (and this was shown earlier by Loray, Pereira, Touzet, \cite{Loray1}, \cite{Loray2}) that the variety of Poisson structures $\mathcal P$ on any smooth cubic threefold $X$ does not depend on the holomorphic structure of $X$ and as a projective variety $\mathcal P \subset {\mathbb P}(H^0(X,{\wedge}^2T_X))\cong {\mathbb P}({\wedge}^2V)$ is isomorphic to the Grassmannian $G(2,V)$ of lines in ${\mathbb P}(V)$. Moreover, $\mathcal P\subset {\mathbb P}({\wedge}^2V)$ is the Pl{\" u}cker embedding.\\

\subsection{Variety of Poisson structures.}

In order to describe equations $[\omega,\omega]=0$ defining the variety of Poisson structures $\mathcal P\subset {\mathbb P}(H^0(X,{\wedge}^2T_X))\cong {\mathbb P}(\mathfrak{so}(5))$, we need to compute the Schouten bracket. Since it is bilinear, it is sufficient to compute numbers 
$$
C_{ijkl}=\frac{1}{2}\cdot [{\epsilon}_{ij},{\epsilon}_{kl}],
$$
where ${\epsilon}_{ij}=Z_jdZ_i-Z_idZ_j, i<j$ are elements of the basis of $H^0(X,{\wedge}^2T_X)$ from Corollary 1.\\

Note that $C_{ijkl}\in H^0(X,{\mathcal O}_X(2))\cong H^0({\mathbb P}^4,{\mathcal O}(2))$.\\

{\bf Lemma 2.} {\it $C_{ijkl}$ is totally antisymmetric with respect to its indices. It is completely determined by the following values:
\begin{gather*}
C_{0123}=\frac{\partial F}{\partial Z_4},\; C_{0124}=-\frac{\partial F}{\partial Z_3},\; C_{0134}=\frac{\partial F}{\partial Z_2}, \; C_{0234}=-\frac{\partial F}{\partial Z_1},\; C_{1234}=\frac{\partial F}{\partial Z_0}.
\end{gather*}
}\\

{\it Proof:} It is enough to work locally assuming, for example, that $Z_0=1$. Let us denote by
$$
X_1=\frac{Z_1}{Z_0},\;\;\;\; X_2=\frac{Z_2}{Z_0},\;\;\;\; X_3=\frac{Z_3}{Z_0},\;\;\;\; X_4=\frac{Z_4}{Z_0}
$$
the remaining affine coordinates and by 
$$
f(X_1,X_2,X_3,X_4)=\frac{F(Z_0,Z_1,Z_2,Z_3,Z_4)}{Z_0^3}
$$
the restriction of the cubic form $F$ defining $X\subset {\mathbb P}^5$. Since $X$ is smooth, without loss of generality we can assume that $\frac{\partial f}{\partial X_4}\neq 0$.\\

Then by Proposition 1 we have 
\begin{multline*}
C_{ijkl}=\frac{1/2}{Vol}\cdot \left( (X_jdX_i-X_idX_j)\wedge (dX_l\wedge dX_k-dX_k\wedge dX_l) + \right. \\
\left. +(dX_j\wedge dX_i-dX_i\wedge dX_j)\wedge (X_ldX_k-X_kdX_l) \right)=\\
=\frac{1}{Vol}\cdot \left( X_i\cdot dX_j\wedge dX_k\wedge dX_l - X_j\cdot dX_i\wedge dX_k\wedge dX_l + X_k\cdot dX_i\wedge dX_j\wedge dX_l - \right.\\
\left. - X_l\cdot dX_i\wedge dX_j\wedge dX_k  \right).
\end{multline*}

The fact that $C_{ijkl}$ is totally antisymmetric is now evident. Hence it is enough to compute $C_{ijkl}$ assuming that $i<j<k<l$.\\

Note that by the adjunction formula we can take
$$
Vol=\frac{dX_1\wedge dX_2 \wedge dX_3}{{\partial f}/{\partial X_4}}.
$$

If $i=0$, then
$$
C_{ijkl}=\frac{dX_j\wedge dX_k \wedge dX_l}{Vol}.
$$

Hence 
\begin{gather*}
C_{0123}=\frac{\partial f}{\partial X_4},\;\; C_{0124}=\frac{\partial X_4}{\partial X_3}\cdot \frac{\partial f}{\partial X_4}=- \frac{\partial f}{\partial X_3},\;\;  C_{0134}=-\frac{\partial X_4}{\partial X_2}\cdot \frac{\partial f}{\partial X_4}= \frac{\partial f}{\partial X_2},\\
C_{0234}=\frac{\partial X_4}{\partial X_1}\cdot \frac{\partial f}{\partial X_4}= -\frac{\partial f}{\partial X_1}.
\end{gather*}

Finally, 
\begin{multline*}
C_{1234}=\frac{\partial f}{\partial X_4}\cdot \frac{1}{dX_1\wedge dX_2 \wedge dX_3}\cdot \left(  X_1\cdot dX_2\wedge dX_3\wedge dX_4 - X_2\cdot dX_1\wedge dX_3\wedge dX_4 +\right. \\
\left. + X_3\cdot dX_1\wedge dX_2\wedge dX_4  - X_4\cdot dX_1\wedge dX_2\wedge dX_3         \right)=\\
=\frac{\partial f}{\partial X_4}\cdot  \left(  X_1\cdot \frac{\partial X_4}{\partial X_1} + X_2\cdot \frac{\partial X_4}{\partial X_2} + X_3\cdot \frac{\partial X_4}{\partial X_3} - X_4 \right)=-\sum_{i=1}^{4} X_i\cdot \frac{\partial f}{\partial X_i}=\frac{\partial F}{\partial Z_0}.
\end{multline*}
{\it QED}\\

Now let $\omega=\sum_{i<j} a_{ij}\cdot {\epsilon}_{ij}$ be a point of ${\mathbb P}(H^0(X,{\wedge}^2T_X))\cong {\mathbb P}(\mathfrak{so}(5))$ with coordinates $(a_{ij})$.\\

Then
$$
[\omega,\omega]=\sum_{i<j} \sum_{k<l} a_{ij}a_{kl}\cdot [{\epsilon}_{ij}, {\epsilon}_{kl}]=4\sum_{i<j<k<l} {\alpha}_{ijkl} \cdot C_{ijkl},
$$
where ${\alpha}_{ijkl}=a_{ij}a_{kl}-a_{ik}a_{jl}+a_{il}a_{jk}$.\\

{\bf Theorem 1 (Loray, Pereira, Touzet, \cite{Loray1}, \cite{Loray2}).} {\it Let $X\subset {\mathbb P}^4$ be a smooth cubic threefold. Then the variety of Poisson structures $\mathcal P\subset {\mathbb P}(\mathfrak{so}(5))$ on $X$ is isomorphic to the Grassmannian $G(2,5)$ of lines in ${\mathbb P}^4$ embedded into ${\mathbb P}(\mathfrak{so}(5))\cong {\mathbb P}({\wedge}^2k^{\oplus 5})$ via the Pl{\" u}cker embedding.}\\

{\it Proof:} It follows from Lemma 2 that 
$$
\frac{1}{4}[\omega,\omega]={\alpha}_{1234}\cdot \frac{\partial F}{\partial Z_0}-{\alpha}_{0234}\cdot \frac{\partial F}{\partial Z_1}+{\alpha}_{0134}\cdot \frac{\partial F}{\partial Z_2}-{\alpha}_{0124}\cdot \frac{\partial F}{\partial Z_3}+{\alpha}_{0123}\cdot \frac{\partial F}{\partial Z_4}.
$$

Since $X\subset {\mathbb P}^4$ is smooth (and $char(\mathbb C)=0$), this expression is $0$ as an element of $H^0(X,{\mathcal O}_X(2))\cong H^0({\mathbb P}^4,{\mathcal O}(2))$ if and only if ${\alpha}_{ijkl}=0$ for any $i<j<k<l$.

Hence $\mathcal P \subset {\mathbb P}(\mathfrak{so}(5))$ is the intersection of the following $5$ quadrics:
$$
{\alpha}_{0123},\; {\alpha}_{0124},\; {\alpha}_{0134},\; {\alpha}_{0234},\; {\alpha}_{1234}.
$$

These are exactly the Pl{\" u}cker quadrics defining $G(2,5)\subset {\mathbb P}^9$. {\it QED}\\

As we mentioned in the introduction, this theorem was proved earlier by Loray, Pereira, Touzet (see \cite{Loray1}) more generally by a different method.\\

\subsection{Poisson cohomology.}

Let us compute Poisson cohomology of $X\subset {\mathbb P}^4$ for any Poisson structure $\omega \in \mathcal P$.\\

According to \cite{Stienon}, Corollary 4.26, Poisson cohomology is equal to the total cohomology of the following double complex ${\Omega}^{\star,\star}$:\\

$$
\begin{CD}
\cdots @. \cdots @. \cdots @. \\
@AAA       @AAA       @AAA  @. \\
{\Omega}^{0,0}(X,{\wedge}^2T_X) @>\bar{\partial}>> {\Omega}^{0,1}(X,{\wedge}^2T_X) @>\bar{\partial}>> {\Omega}^{0,2}(X,{\wedge}^2T_X) @>\bar{\partial}>> \cdots \\
@AAd_{\omega}A  @AAd_{\omega}A  @AAd_{\omega}A  @. \\
{\Omega}^{0,0}(X,T_X) @>\bar{\partial}>> {\Omega}^{0,1}(X,T_X) @>\bar{\partial}>> {\Omega}^{0,2}(X,T_X) @>\bar{\partial}>> \cdots \\
@AAd_{\omega}A  @AAd_{\omega}A  @AAd_{\omega}A  @. \\
{\Omega}^{0,0}(X,{\mathcal O}_X) @>\bar{\partial}>> {\Omega}^{0,1}(X,{\mathcal O}_X) @>\bar{\partial}>> {\Omega}^{0,2}(X,{\mathcal O}_X) @>\bar{\partial}>> \cdots
\end{CD}
$$

Here the vertical maps $d_{\omega}\colon {\wedge}^iT_X\rightarrow {\wedge}^{i+1}T_X$ are given by the Schouten bracket with $\omega$:
$$
d_{\omega}(\nu)=\frac{1}{2}[\omega,\nu].
$$
\par
Hence Poisson cohomology $H^{\star}_{Poisson}(X,\omega)$ is computed by the spectral sequence of this double complex:
$$
E_2^{p,q}= {^{v}H}^{p} {^{h}H}^{q}({\Omega}^{\star,\star}) \Rightarrow H^{p+q}_{Poisson}(X,\omega).
$$

The cohomology groups of horizontal complexes in the double complex ${\Omega}^{\star,\star}$ are exactly the cohomology groups of sheaves of multivector fields on $X$.\\

Hence from Lemma 1 it follows that $^{h}H^{q}({\Omega}^{p,\star})$ has the following shape:\\
\begin{center}
\begin{tabular}{ccc}
0 & 0 & 0 \\
$H^0(X,{\mathcal O}_X(2))$ & 0 & 0 \\
$\uparrow d_{\omega}$ & & \\
$H^0(X,{\wedge}^2T_X)$ & 0 & 0 \\
0 & $H^1(X,T_X)$ & 0 \\
$H^0(X,{\mathcal O}_X)$ & 0 & 0 \\
\end{tabular}
\end{center}

The matrix $C_{\omega}$ of the linear map
$$
d_{\omega}\colon \mathfrak{so}(5)\cong H^0(X,{\wedge}^2T_X)\rightarrow H^0(X,{\mathcal O}_X(2))\cong H^0({\mathbb P}^4,{\mathcal O}(2))
$$
can be computed explicitly. This is done in the next Proposition 2.\\

{\bf Proposition 2.} {\it Let $\omega=\sum_{i<j}a_{ij}\cdot {\epsilon}_{ij}\in \mathcal P$. Then the images of the basis elements ${\epsilon}_{kl}\in \mathfrak{so}(5)$ under $d_{\omega}$ are as follows:
\begin{gather*}
d_{\omega}({\epsilon}_{01})=a_{23}\frac{\partial F}{\partial Z_4}-a_{24}\frac{\partial F}{\partial Z_3}+a_{34}\frac{\partial F}{\partial Z_2},\;  d_{\omega}({\epsilon}_{02})=-a_{13}\frac{\partial F}{\partial Z_4}+a_{14}\frac{\partial F}{\partial Z_3}-a_{34}\frac{\partial F}{\partial Z_1}, \\
d_{\omega}({\epsilon}_{03})=a_{12}\frac{\partial F}{\partial Z_4}-a_{14}\frac{\partial F}{\partial Z_2}+a_{24}\frac{\partial F}{\partial Z_1},\;  d_{\omega}({\epsilon}_{04})=-a_{12}\frac{\partial F}{\partial Z_3}+a_{13}\frac{\partial F}{\partial Z_2}-a_{23}\frac{\partial F}{\partial Z_1}, \\
d_{\omega}({\epsilon}_{12})=a_{03}\frac{\partial F}{\partial Z_4}-a_{04}\frac{\partial F}{\partial Z_3}+a_{34}\frac{\partial F}{\partial Z_0},\;  d_{\omega}({\epsilon}_{13})=-a_{02}\frac{\partial F}{\partial Z_4}+a_{04}\frac{\partial F}{\partial Z_2}-a_{24}\frac{\partial F}{\partial Z_0}, \\
d_{\omega}({\epsilon}_{14})=a_{02}\frac{\partial F}{\partial Z_3}-a_{03}\frac{\partial F}{\partial Z_2}+a_{23}\frac{\partial F}{\partial Z_0},\;  d_{\omega}({\epsilon}_{23})=a_{01}\frac{\partial F}{\partial Z_4}-a_{04}\frac{\partial F}{\partial Z_1}+a_{14}\frac{\partial F}{\partial Z_0}, \\
d_{\omega}({\epsilon}_{24})=-a_{01}\frac{\partial F}{\partial Z_3}+a_{03}\frac{\partial F}{\partial Z_1}-a_{13}\frac{\partial F}{\partial Z_0},\;  d_{\omega}({\epsilon}_{34})=a_{01}\frac{\partial F}{\partial Z_2}-a_{02}\frac{\partial F}{\partial Z_1}+a_{12}\frac{\partial F}{\partial Z_0}.
\end{gather*}
}\\

{\it Proof:} $[\omega,{\epsilon}_{kl}]=\sum_{i<j}a_{ij}\cdot [{\epsilon}_{ij},{\epsilon}_{kl}]=2\cdot \sum_{i<j}a_{ij}\cdot C_{ijkl}$.\\

In particular,
\begin{gather*}
\frac{1}{2}[\omega,{\epsilon}_{01}]=a_{23}C_{0123}+a_{24}C_{0124}+a_{34}C_{0134},\;  \frac{1}{2}[\omega,{\epsilon}_{02}]=-a_{13}C_{0123}-a_{14}C_{0124}+a_{34}C_{0234}, \\
\frac{1}{2}[\omega,{\epsilon}_{03}]=a_{12}C_{0123}-a_{14}C_{0134}-a_{24}C_{0234},\;   \frac{1}{2}[\omega,{\epsilon}_{04}]=a_{12}C_{0124}+a_{13}C_{0134}+a_{23}C_{0234}, \\
\frac{1}{2}[\omega,{\epsilon}_{12}]=a_{03}C_{0123}+a_{04}C_{0124}+a_{34}C_{1234},\;   \frac{1}{2}[\omega,{\epsilon}_{13}]=-a_{02}C_{0123}+a_{04}C_{0134}-a_{24}C_{1234}, \\
\frac{1}{2}[\omega,{\epsilon}_{14}]=-a_{02}C_{0124}-a_{03}C_{0134}+a_{23}C_{1234},\;   \frac{1}{2}[\omega,{\epsilon}_{23}]=a_{01}C_{0123}+a_{04}C_{0234}+a_{14}C_{1234}, \\
\frac{1}{2}[\omega,{\epsilon}_{24}]=a_{01}C_{0124}-a_{03}C_{0234}-a_{13}C_{1234},\;   \frac{1}{2}[\omega,{\epsilon}_{34}]=a_{01}C_{0134}+a_{02}C_{0234}+a_{12}C_{1234}.
\end{gather*}

Now one uses Lemma 2. {\it QED}\\

Hence the second sheet $E_2^{\star,\star}$ of the Laurent-Gengoux-Sti{\' e}non-Xu spectral sequence has the following shape:\\

\begin{center}
\begin{tabular}{ccc}
0 & 0 & 0 \\
$coker(d_{\omega})$ & 0 & 0 \\
$ker(d_{\omega})$ & 0 & 0 \\
0 & $H^1(X,T_X)$ & 0 \\
$k$ & 0 & 0 \\
\end{tabular}
\end{center}

Since all the differentials vanish, we have
$$
E_{\infty}^{p,q}=E_2^{p,q}.
$$

This proves the following result.\\

{\bf Theorem 2.} {\it Let $X\subset {\mathbb P}^4$ be a smooth cubic threefold and $\omega\in \mathcal P\subset {\mathbb P}(H^0(X,{\wedge}^2T_X))$ a Poisson structure on $X$. Then 
\begin{itemize}
\item $H^0_{Poisson}(X,\omega)\cong H^0(X,{\mathcal O}_X)=k$,
\item $H^1_{Poisson}(X,\omega)=0$,
\item $H^2_{Poisson}(X,\omega)$ is an extension of $ker(d_{\omega})$ by $H^1(X,T_X)$,
\item $H^3_{Poisson}(X,\omega)\cong coker(d_{\omega})$.
\end{itemize}}

In particular, 
$$
dim(H^2_{Poisson}(X,\omega))=dim(ker(d_{\omega})) + dim(H^1(X,T_X))=10-rk(C_{\omega})+10=20-rk(C_{\omega})
$$
and 
$$
dim(H^3_{Poisson}(X,\omega))=15-rk(C_{\omega}).
$$

\section{Del Pezzo quintic threefold.}

Let $X$ be the (smooth) del Pezzo quintic threefold. For the standard properties of $X$ (as well as for Fano varieties in general) we refer the reader to \cite{Fano} or papers of Iskovskikh, Mukai,...\\

The most important for us description of $X$ is given as a general codimension $3$ linear section of the Grassmannian $G(2,5)\subset {\mathbb P}^9$ of lines in ${\mathbb P}^4$ in its Pl{\" u}cker embedding.\\

{\bf Remark.} In what follows we use essentially the same approach (via cohomology and spaces of global sections of vector bundles) as in Section 2 in order to work with Poisson structures on $X$. For this the description of $X$ as a linear section of the Grassmannian is essential. Alternatively, one can use the fact that $X$ is a compactification of ${\mathbb C}^3$. This allows one to work on ${\mathbb C}^3$ instead of $X$ (provided that one keeps track of which (bi)vector fields extend to the whole of $X$). Indeed, using, for example, \cite{Kimura}, section 11 it is straightforward to describe $X$ in terms of equations in ${\mathbb P}^6$ and find explicitly a scroll (with one double line) on $X$, whose complement is exactly ${\mathbb C}^3$.\\

Let us denote by $Z_0,Z_1,Z_2,Z_3,Z_4,Z_5,Z_6,Z_7,Z_8,Z_9$ the homogeneous coordinates on ${\mathbb P}^9$. Then $G(2,5)\subset {\mathbb P}^9$ is given by the intersection of the following quadrics (see \cite{Mukai}, for example):
\begin{gather*}
p_1=Z_0Z_7-Z_1Z_5+Z_2Z_4,\;\; p_2=Z_0Z_8-Z_1Z_6+Z_3Z_4,\;\; p_3=Z_0Z_9-Z_2Z_6+Z_3Z_5,\\
p_4=Z_1Z_9-Z_2Z_8+Z_3Z_7,\;\;\;\;\;\; p_5=Z_4Z_9-Z_8Z_5+Z_6Z_7.
\end{gather*}

We will use the following hyperplanes in order to obtain $X$ in the intersection with the Grassmannian (see \cite{Kimura}, section 11):
\begin{gather*}
{\lambda}_1=Z_0+Z_7,\;\; {\lambda}_2=Z_4+Z_9,\;\; {\lambda}_3=Z_1+Z_6.
\end{gather*}

(The further intersection with one of the coordinate hyperplanes gives the scroll mentioned in the Remark above, whose complement is ${\mathbb C}^3$.)\\

\subsection{Cohomology computations.}

Let $G=G(2,5)\subset {\mathbb P}^9$ and ${\Pi}_i\subset {\mathbb P}^9, i=1,2,3$ be the hyperplanes (given by the equations ${\lambda}_i=0$ as above). Let $X_0=G$, $X_1=G\cap {\Pi}_1$, $X_2=G\cap {\Pi}_1 \cap {\Pi}_2$, $X_3=G\cap {\Pi}_1 \cap {\Pi}_2 \cap {\Pi}_3=X$. Then we have a chain of embeddings
$$
X=X_3\subset X_2\subset X_1 \subset X_0=G\subset {\mathbb P}^9.
$$

As in the case of cubic threefolds, we need to show vanishing of higher cohomology groups and describe the spaces of global sections of sheaves of multivector fields on $X$.\\

Let $S$ and $Q$ denote the tautological (of rank $2$) and the quotient (of rank $3$) vector bundles on $G=G(2,5)$.\\

We will use the exact sequence
\begin{equation}
\label{TES1}
0\rightarrow {\mathcal O}_{X_i}(-1)\rightarrow {\mathcal O}_{X_i} \rightarrow {\mathcal O}_{X_{i+1}} \rightarrow 0 
\end{equation}

on $X_i$, the tautological exact sequence

\begin{equation}
\label{TES2}
0\rightarrow S \rightarrow {\mathcal O}_{G}^{\oplus 5} \rightarrow Q \rightarrow 0 
\end{equation}

on $G$ and the conormal bundle sequence 

\begin{equation}
\label{TES1pp}
0\rightarrow N_{X/G}^{\vee} \rightarrow {\Omega}_{G}^{1} {\mid}_X \rightarrow {\Omega}_{X}^{1} \rightarrow 0 
\end{equation}
on X.

In order to compute cohomology of various vector bundles on $G$ we will use general formulas from \cite{Kapranov} (section 3) and \cite{Fonarev} (section 2.2). In the notation of \cite{Kapranov} $S^{\bot}=Q^{\vee}$.\\

Note that $T_G\cong S^{\vee}\otimes Q$, ${\Omega}^1_{G}\cong S\otimes Q^{\vee}$, $N_{G/{\mathbb P}^9}\cong Q^{\vee} (2)$ (see \cite{Manivel}, section 5.4) and $N_{X/G}\cong \oplus_{i=1}^3 {\mathcal O}_X(1)$.\\

In the next Lemma we collect some vanishing results for cohomology of certain vector bundles on the Grassmannian $G=G(2,5)$. They all are consequences of the Borel-Weil-Bott theorem (via the general formulas of Fonarev \cite{Fonarev} and Kapranov \cite{Kapranov}) and are well-known.\\

{\bf Lemma 3.} {\it Let $G=G(2,5)$. Then 
\begin{itemize}
\item[(a)] $H^i(G,{\mathcal O}(-d))=0$ for any $i\geq 1$, $d\leq 3$,
\item[(b)] $H^i(G,{\Omega}^1_G (2-k))=0$ for any $i\geq 1$, $k=0,1,3,4$, $H^i(G,{\Omega}^1_G )=0$ for $i\neq 1$ and $H^1(G,{\Omega}^1_G)\cong k$,
\item[(c)] $H^0(G,{\Omega}^1_G (d))=0$ for $d=-1,0,1$,
\item[(d)] $H^0(G,{\Omega}^1_G (2))\cong k^{\oplus 45}$,
\item[(e)] $H^1(G,{\mathcal O}(d))=0$ for any $d\in \mathbb Z$,
\item[(f)] $H^0(G,{\mathcal O}(2))\cong k^{\oplus 50}$.
\end{itemize}}

{\it Proof:} In the notation of \cite{Fonarev} and \cite{Kapranov} ${\mathcal O}(-d)={\Sigma}^{-d,-d}S^{\vee}$ and ${\Omega}^1_{G}(d)={\Sigma}^{d,d-1}S^{\vee}\otimes {\Sigma}^{1,0,0}S^{\bot}$.\\

(a) follows, for example, from \cite{Kapranov}, Lemma 3.2 and (b), (c), (e) follow from \cite{Fonarev}, section 2.2.\\

(d) According to \cite{Fonarev}, the dual of $H^0(G, {\Omega}^1_G (2))$ admits an irreducible representation of $GL(5)$ corresponding to the Young diagram 
\begin{Young}
&\cr
\cr
\cr
\end{Young}. By the hook formula, its dimension is $45$.\\

(f) According to \cite{Fonarev}, section 2.2, the dual of $H^0(G, {\mathcal O}(2))$ admits an irreducible representation of $GL(5)$ corresponding to the Young diagram 
\begin{Young}
&\cr
&\cr
\end{Young}. By the hook formula, its dimension is $50$. {\it QED}\\

The following observation will be also employed.\\

{\bf Remark.} Let $\mathcal E$ be a vector bundle on $G$. Then the vanishing of $H^i(X,\mathcal E)$ follows from the vanishing of $H^{i+k}(G,{\mathcal E}(-k))$ for any $k=0,1,2,3$. Indeed, one multiplies~(\ref{TES1}) by $\mathcal E$, takes cohomology and applies induction.\\

{\bf Lemma 4.} {\it Let $X$ be the (smooth) del Pezzo quintic threefold. Then 
\begin{itemize}
\item[(a)] $H^0(X,{\mathcal O}_X)=k$, $H^i(X,{\mathcal O}_X)=H^i(X,{\mathcal O}_X(1))=0$ for any $i\geq 1$,
\item[(b)] $H^0(X,T_X)\cong k^{\oplus 3}$, $H^i(X,T_X)=0$ for any $i\geq 1$,
\item[(c)] $H^0(X,{\wedge}^2T_X)\cong k^{\oplus 21}$, $H^i(X,{\wedge}^2T_X)=0$ for any $i\geq 1$,
\item[(d)] $H^0(X,-K_X)\cong k^{\oplus 23}$, $H^i(X,-K_X)=0$ for any $i\geq 1$.
\end{itemize}}

{\it Proof:} (a) This follows from the Remark above and Lemma 3(a).\\

(b) $H^i(X,T_X)\cong H^i(X,{\Omega}^2_X(2))=0$ for $i\geq 2$ by the Kodaira Vanishing theorem. $H^0(X,T_X)\cong k^{\oplus 3}$ is well-known (see \cite{Mukai2}, for example).\\

(c) Note that ${\wedge}^2T_X\cong {\Omega}^1_X(2)$.\\

Multiplying~(\ref{TES1pp}) by ${\mathcal O}(2)$ and taking cohomology, we obtain the short exact sequence
\begin{equation}
\label{TES5}
0\rightarrow H^0(X,{\mathcal O}(1))^{\oplus 3} \rightarrow H^0(X,{\Omega}_{G}^{1}(2)) \rightarrow H^0(X,{\Omega}_{X}^{1}(2)) \rightarrow 0 
\end{equation}

and isomorphisms
$$
H^i(X,{\Omega}_{G}^{1}(2))\cong H^i(X,{\Omega}_{X}^{1}(2))\;\; \mbox{for any}\;\; i\geq 1.
$$
By the Remark above the vanishing of $H^i(X,{\Omega}_{G}^{1}(2))$, $i\geq 1$ follows from the vanishing of $H^{i+k}(G,{\Omega}_{G}^{1}(2-k))$ for any $k=0,1,2,3$ and any $i\geq 1$. This follows from Lemma 3(b).\\

Moreover, by~(\ref{TES5})
$$
dim(H^0(X,{\Omega}_{X}^{1}(2)))=dim(H^0(X,{\Omega}_{G}^{1}(2)))-3dim(H^0(X,{\mathcal O}(1))).
$$
It follows from part (a) that $dim(H^0(X,{\mathcal O}(1)))=dim(H^0(G,{\mathcal O}(1)))-3=7$. Let us compute $dim(H^0(X,{\Omega}_{G}^{1}(2)))$.\\

Taking cohomology of the short exact sequences
\begin{gather*}
0\rightarrow {\Omega}_{G}^{1}(1) {\mid}_{X_2} \rightarrow {\Omega}_{G}^{1}(2) {\mid}_{X_2} \rightarrow {\Omega}_{G}^{1}(2) {\mid}_{X} \rightarrow 0,\\ 
0\rightarrow {\Omega}_{G}^{1}(1) {\mid}_{X_1} \rightarrow {\Omega}_{G}^{1}(2) {\mid}_{X_1} \rightarrow {\Omega}_{G}^{1}(2) {\mid}_{X_2} \rightarrow 0,\\ 
0\rightarrow {\Omega}_{G}^{1}(1) \rightarrow {\Omega}_{G}^{1}(2) \rightarrow {\Omega}_{G}^{1}(2) {\mid}_{X_1} \rightarrow 0,\\ 
0\rightarrow {\Omega}_{G}^{1} {\mid}_{X_1} \rightarrow {\Omega}_{G}^{1}(1) {\mid}_{X_1} \rightarrow {\Omega}_{G}^{1}(1) {\mid}_{X_2} \rightarrow 0,\\
0\rightarrow {\Omega}_{G}^{1} \rightarrow {\Omega}_{G}^{1}(1) \rightarrow {\Omega}_{G}^{1}(1) {\mid}_{X_1} \rightarrow 0,\\ 
0\rightarrow {\Omega}_{G}^{1}(-1) \rightarrow {\Omega}_{G}^{1} \rightarrow {\Omega}_{G}^{1} {\mid}_{X_1} \rightarrow 0, 
\end{gather*}

and applying the vanishing observations from Lemma 3, one obtains the following example sequences:
\begin{gather*}
0\rightarrow H^0(X_2,{\Omega}_{G}^{1}(1)) \rightarrow H^0(X_2,{\Omega}_{G}^{1}(2)) \rightarrow H^0(X,{\Omega}_{G}^{1}(2)) \rightarrow 0,\\ 
0\rightarrow H^0(X_1,{\Omega}_{G}^{1}(1)) \rightarrow H^0(X_1,{\Omega}_{G}^{1}(2)) \rightarrow H^0(X_2,{\Omega}_{G}^{1}(2)) \rightarrow 0,\\ 
0\rightarrow H^0( G, {\Omega}_{G}^{1}(2)) \rightarrow H^0(X_1,{\Omega}_{G}^{1}(2)) \rightarrow 0,\\ 
0\rightarrow H^0(X_1,{\Omega}_{G}^{1}) \rightarrow H^0(X_1,{\Omega}_{G}^{1}(1)) \rightarrow H^0(X_2,{\Omega}_{G}^{1}(1)) \rightarrow H^1(X_1,{\Omega}_{G}^{1}) \rightarrow 0,\\
0\rightarrow H^0(X_1,{\Omega}_{G}^{1}(1))\rightarrow H^1(G,{\Omega}_{G}^{1}) \rightarrow 0,\\ 
0\rightarrow H^0(X_1,{\Omega}_{G}^{1}) \rightarrow 0,\\
0\rightarrow H^1(G,{\Omega}_{G}^{1})\rightarrow H^1(X_1,{\Omega}_{G}^{1}) \rightarrow 0.
\end{gather*}

This implies that 
$$
dim(H^0(X,{\Omega}_{G}^{1}(2)))=dim(H^0(X_1,{\Omega}_{G}^{1}(2)))-3 dim (H^0(X_1,{\Omega}_{G}^{1}(1)))=dim (H^0(G,{\Omega}_{G}^{1}(2)))-3=42.
$$

Hence $dim(H^0(X,{\Omega}_{X}^{1}(2)))=42-3\cdot 7=21$.\\

(d) Since $\mathcal O (-K_X)\cong {\mathcal O}(2)$ on $X$, part (d) follows from the Remark above and Lemma 3.\\

In order to compute $dim(H^0(X,{\mathcal O}(2)))$, one can use the same method as in part (c), i.e. the induction along the chain $X=X_3\subset X_2\subset X_1\subset X_0=G$. Then using Lemma 3 one obtains that 
$$
dim(H^0(X,{\mathcal O}(2)))=dim(H^0(G,{\mathcal O}(2)))-3\cdot dim(H^0(G,{\mathcal O}(1)))+3=50-3\cdot 10+3=23.
$$
{\it QED}\\

Now let us give descriptions of the spaces of global sections $H^0(X,T_X)$, $H^0(X,{\wedge}^2T_X)$ and $H^0(X,-K_X)$.\\

Let $X$ be the del Pezzo quintic threefold as above (i.e. the intersection of quadrics $p_1,p_2,p_3,p_4,p_5$ and hyperplanes ${\lambda}_1,{\lambda}_2,{\lambda}_3$ in ${\mathbb P}^9$).\\

{\bf Lemma 5.} {\it As a basis of $H^0(X,T_X)\cong \mathfrak{so}(2)$ one can take the restrictions to $X\subset G\subset {\mathbb P}^9$ of the following vector fields on ${\mathbb P}^9$:
\begin{itemize}
\item $v_1=2Z_1\cdot (\frac{\partial}{\partial Z_1}-\frac{\partial}{\partial Z_6})-Z_2\cdot \frac{\partial}{\partial Z_2}+3Z_3\cdot \frac{\partial}{\partial Z_3}+Z_4\cdot (\frac{\partial}{\partial Z_4}-\frac{\partial}{\partial Z_9})-2Z_5\cdot \frac{\partial}{\partial Z_5}+4Z_8\cdot \frac{\partial}{\partial Z_8}$,
\item $v_2=Z_2\cdot (\frac{\partial}{\partial Z_0}-\frac{\partial}{\partial Z_7})+3Z_4\cdot (\frac{\partial}{\partial Z_1}-\frac{\partial}{\partial Z_6})+3Z_5\cdot \frac{\partial}{\partial Z_2}-5Z_1\cdot \frac{\partial}{\partial Z_3}+2Z_0\cdot (\frac{\partial}{\partial Z_4}-\frac{\partial}{\partial Z_9})-Z_3\cdot \frac{\partial}{\partial Z_8}$,
\item $v_3=-3Z_4\cdot (\frac{\partial}{\partial Z_0}-\frac{\partial}{\partial Z_7})+Z_3\cdot(\frac{\partial}{\partial Z_1}-\frac{\partial}{\partial Z_6})-5Z_0\cdot \frac{\partial}{\partial Z_2}+3Z_8\cdot \frac{\partial}{\partial Z_3}-2Z_1\cdot (\frac{\partial}{\partial Z_4}-\frac{\partial}{\partial Z_9})-Z_2\cdot \frac{\partial}{\partial Z_5}$.
\end{itemize}}

{\it Proof:} The fact that $v_i$ restrict to vector fields on $X$ follows from the vanishing
$$
v_i(dp_j)=0, \;\;\; j=1,2,3,4,5\;\;\mbox{and }\;\; v_i(d{\lambda}_k)=0, \;\;\; k=1,2,3.
$$

The fact that $v_1,v_2,v_3$ are linearly independent can be checked locally on $X$ (in an open set $Z_8=1$, for example - see Lemma 8 below). {\it QED}\\

{\bf Lemma 6.} {\it One can identify $H^0(X,{\wedge}^2T_X)$ with the Lie algebra (viewed merely as a vector space) $\mathfrak{so}(7)$ of skew-symmetric $7\times 7$ matrices.\\

(Or with ${\wedge}^2V$, where $V=H^0(X,{\mathcal O}(1))$.)\\

As a basis of $H^0(X,{\wedge}^2T_X)\cong H^0(X,{\Omega}^1_X(2))$ one can take the restrictions to $X\subset G\subset {\mathbb P}^9$ of the forms ${\epsilon}_{ij}=Z_jdZ_i-Z_idZ_j$ on ${\mathbb P}^9$, where $0\leq i < j \leq 5$ or $0\leq i \leq 5$, $j=8$.}\\

{\it Proof:} The fact that the restrictions of ${\epsilon}_{ij}$ to $X$ give rise to $21$ linearly independent global sections of ${\Omega}^1_X(2)$ can be checked locally on $X$ (in an open set $Z_8=1$, for example - see Lemma 8 below).

Since $dim(H^0(X,{\wedge}^2T_X))=21$ by Lemma 4, the result follows. {\it QED}\\

From now on we will always assume that the range of indices $i<j$ of ${\epsilon}_{ij}$ is the same as in Lemma 6.\\

{\bf Lemma 7.} {\it As a basis of $H^0(X,-K_X)\cong H^0(X,{\mathcal O}_X(2))$ one can take the restrictions to $X\subset G\subset {\mathbb P}^9$ of the quadratic forms ${z}_{ij}=Z_iZ_j$ on ${\mathbb P}^9$, where $0\leq i \leq j \leq 9$, $i,j \neq 6,7,9$ and $(ij)\neq (04),(14),(24),(34),(44)$.}\\

{\it Proof:} Since $dim(H^0(X,{\mathcal O}_X(2)))=23$ and we are given precisely $23$ elements $z_{ij}$, it is enough to check locally that they are linearly independent. {\it QED}\\

\subsection{Matrix elements of Schouten brackets.}

Given an element $\omega\in H^0(X,{\wedge}^2T_X)$, we will need to know the matrices of the linear maps
$$
{\alpha}_{\omega}\colon H^0(X,T_X)\rightarrow H^0(X,{\wedge}^2T_X)\cong H^0(X,{\Omega}^1_X(2)), \;\;\; \nu\mapsto [\nu ,\omega]
$$
and 
$$
{\beta}_{\omega}\colon H^0(X,{\wedge}^2T_X)\rightarrow H^0(X,-K_X)\cong H^0(X,{\mathcal O}(2)), \;\;\; \nu\mapsto [\nu ,\omega]
$$

Let us denote these matrices by $A_{\omega}$ and $B_{\omega}$ respectively. In this subsection we compute them relative to the bases in $H^0(X,{\wedge}^jT_X)$ introduced in Lemma 5, Lemma 6 and Lemma 7.\\

In order to do this, it is sufficient to compute the images ${\alpha}_{\omega}(v_i)$, $i=1,2,3$ and ${\beta}_{\omega}({\epsilon}_{ij})$, $i<j$.\\

If $\omega=\sum_{i<j}a_{ij}\cdot {\epsilon}_{ij}$, then 
$$
{\alpha}_{\omega}(v_i)=\sum_{j<k}a_{jk}\cdot [v_i,{\epsilon}_{jk}]\;\; \mbox{and}\;\; {\beta}_{\omega}({\epsilon}_{ij})=\sum_{k<l}a_{kl}\cdot [{\epsilon}_{ij},{\epsilon}_{kl}].
$$

Hence $A_{\omega}$ and $B_{\omega}$ are determined by $A_{ijk}=[v_i,{\epsilon}_{jk}]$ and $B_{ijkl}=\frac{1}{2}[{\epsilon}_{ij},{\epsilon}_{kl}]$. Let us compute them.\\

{\bf Lemma 8.} {\it
\begin{gather*}
A_{101}=0,\; A_{102}=-3{\epsilon}_{02},\; A_{103}={\epsilon}_{03},\; A_{104}=-{\epsilon}_{04},\; A_{105}=-4{\epsilon}_{05},\; A_{108}=2{\epsilon}_{08},\\
A_{112}=-{\epsilon}_{12},\; A_{113}=3{\epsilon}_{13},\; A_{114}={\epsilon}_{14},\; A_{115}=-2{\epsilon}_{15},\; A_{118}=4{\epsilon}_{18},\; A_{123}=0,\\
A_{124}=-2{\epsilon}_{24},\; A_{125}=-5{\epsilon}_{25},\; A_{128}={\epsilon}_{28},\; A_{134}=2{\epsilon}_{34},\; A_{135}=-{\epsilon}_{35},\\
A_{138}=5{\epsilon}_{38},\; A_{145}=-3{\epsilon}_{45},\; A_{148}=3{\epsilon}_{48},\; A_{158}=0,\; A_{201}=3{\epsilon}_{04}-{\epsilon}_{12},\\
A_{202}=3{\epsilon}_{05},\; A_{203}=-5{\epsilon}_{01}+{\epsilon}_{23},\; A_{204}={\epsilon}_{24},\; A_{205}={\epsilon}_{25},\; A_{208}=-{\epsilon}_{03}+{\epsilon}_{28},\\
A_{212}=3{\epsilon}_{15}-3{\epsilon}_{24},\; A_{213}=-3{\epsilon}_{34},\; A_{214}=-2{\epsilon}_{01},\; A_{215}=3{\epsilon}_{45},\\
A_{218}=-{\epsilon}_{13}+3{\epsilon}_{48},\; A_{223}=5{\epsilon}_{12}-3{\epsilon}_{35},\; A_{224}=-2{\epsilon}_{02}-3{\epsilon}_{45},\; A_{225}=0,\\
A_{228}=-{\epsilon}_{23}+3{\epsilon}_{58},\; A_{234}=-3{\epsilon}_{03}-5{\epsilon}_{14},\; A_{235}=-5{\epsilon}_{15},\; A_{238}=-5{\epsilon}_{18},\\
A_{245}=2{\epsilon}_{05},\; A_{248}=3{\epsilon}_{08}+{\epsilon}_{34},\; A_{258}={\epsilon}_{35},\; A_{301}={\epsilon}_{03}+3{\epsilon}_{14},\; A_{302}=3{\epsilon}_{24},\\
A_{303}=3{\epsilon}_{08}+3{\epsilon}_{34},\; A_{304}=-2{\epsilon}_{01},\; A_{305}=-{\epsilon}_{02}-3{\epsilon}_{45},\; A_{308}=-3{\epsilon}_{48},\\
A_{312}=5{\epsilon}_{01}-{\epsilon}_{23},\; A_{313}=3{\epsilon}_{18},\; A_{314}={\epsilon}_{34},\; A_{315}=-{\epsilon}_{12}+{\epsilon}_{35},\; A_{318}={\epsilon}_{38},\\
A_{323}=-5{\epsilon}_{03}+3{\epsilon}_{28},\; A_{324}=-5{\epsilon}_{04}+2{\epsilon}_{12},\; A_{325}=-5{\epsilon}_{05},\; A_{328}=-5{\epsilon}_{08},\\
A_{334}=2{\epsilon}_{13}-3{\epsilon}_{48},\; A_{335}={\epsilon}_{23}-3{\epsilon}_{58}, A_{338}=0,\; A_{345}=-2{\epsilon}_{15}+{\epsilon}_{24},\\
A_{348}=-2{\epsilon}_{18},\; A_{358}=-{\epsilon}_{28}.
\end{gather*}}

{\it Proof:} It is enough to work locally. Without loss of generality we can assume that $Z_8=1$. Let us denote by $x_i=\frac{Z_i}{Z_8}, 0\leq i\leq 9, i\neq 8$ the affine coordinates on this open subset of ${\mathbb P}^9$.\\

Over this open set $X$ is isomorphic to the affine space ${\mathbb C}^3$ with coordinates $x_1,x_3,x_4$ and we have on $X$ the following relations:
\begin{gather*}
x_0=-x_1^2-x_3x_4,\;\; x_2=-x_1x_4+x_3x_1^2+x_3^2x_4,\;\; x_5=-x_1^3-x_4^2-x_1x_3x_4,\\
x_9=-x_4,\;\;\;\;\; x_7=-x_0,\;\;\;\;\; x_6=-x_1.
\end{gather*}

The restrictions of the vector fields $v_i$, $i=1,2,3$ from Lemma 5 to this open subset have the following form:
\begin{itemize}
\item $v_1=-2x_1\frac{\partial}{\partial x_1}-x_3\frac{\partial}{\partial x_3}-3x_4\frac{\partial}{\partial x_4},$
\item $v_2=(x_1x_3+3x_4)\frac{\partial}{\partial x_1}+(x_3^2-5x_1)\frac{\partial}{\partial x_3}-(2x_1^2+x_3x_4)\frac{\partial}{\partial x_4},$
\item $v_3=x_3\frac{\partial}{\partial x_1}+3\frac{\partial}{\partial x_3}-2x_1\frac{\partial}{\partial x_4}.$
\end{itemize}

The restrictions of 1-forms ${\epsilon}_{ij}=Z_jdZ_i-Z_idZ_j$, $i<j$ from Lemma 6 to this open subset are just $x_jdx_i-x_idx_j$. One obtains:
\begin{itemize}
\item ${\epsilon}_{01}=(x_3x_4-x_1^2)dx_1-x_1x_4dx_3-x_1x_3dx_4,$
\item ${\epsilon}_{02}=(x_1^2x_4-x_3x_4^2)dx_1+(x_1^4+x_3^2x_4^2+2x_1^2x_3x_4+x_1x_4^2)dx_3-x_1^3dx_4,$
\item ${\epsilon}_{03}=-2x_1x_3dx_1+x_1^2dx_3-x_3^2dx_4,$
\item ${\epsilon}_{04}=-2x_1x_4dx_1-x_4^2dx_3+x_1^2dx_4,$
\item ${\epsilon}_{05}=(2x_1x_4^2-2x_1^2x_3x_4-x_1^4-x_3^2x_4^2)dx_1+x_4^3dx_3-(x_3x_4^2+2x_1^2x_4)dx_4,$
\item ${\epsilon}_{08}=-2x_1dx_1-x_4dx_3-x_3dx_4,$
\item ${\epsilon}_{12}=(x_3^2x_4-x_1^2x_3)dx_1-(x_1^3+2x_1x_3x_4)dx_3+(x_1^2-x_1x_3^2)dx_4,$
\item ${\epsilon}_{13}=x_3dx_1-x_1dx_3,$
\item ${\epsilon}_{14}=x_4dx_1-x_1dx_4,$
\item ${\epsilon}_{15}=(2x_1^3-x_4^2)dx_1+x_1^2x_4dx_3+(x_1^2x_3+2x_1x_4)dx_4,$
\item ${\epsilon}_{18}=dx_1,$
\item ${\epsilon}_{23}=(2x_1x_3^2-x_3x_4)dx_1+(x_3^2x_4+x_1x_4)dx_3+(x_3^3-x_1x_3)dx_4,$
\item ${\epsilon}_{24}=(2x_1x_3x_4-x_4^2)dx_1+(x_1^2x_4+2x_3x_4^2)dx_3-x_1^2x_3dx_4,$
\item ${\epsilon}_{25}=(x_1^4x_3-2x_1x_3x_4^2+x_4^3-2x_1^3x_4+2x_1^2x_3^2x_4+x_3^3x_4^2)dx_1-(x_1^5+2x_1^2x_4^2+2x_1^3x_3x_4+2x_3x_4^3+x_1x_3^2x_4^2)dx_3+(x_1^4+2x_1^2x_3x_4-x_1x_4^2+x_3^2x_4^2)dx_4,$
\item ${\epsilon}_{28}=(2x_1x_3-x_4)dx_1+(x_1^2+2x_3x_4)dx_3+(x_3^2-x_1)dx_4,$
\item ${\epsilon}_{34}=x_4dx_3-x_3dx_4,$
\item ${\epsilon}_{35}=(3x_1^2x_3+x_3^2x_4)dx_1-(x_1^3+x_4^2)dx_3+(2x_3x_4+x_1x_3^2)dx_4,$
\item ${\epsilon}_{38}=dx_3,$
\item ${\epsilon}_{45}=(3x_1^2x_4+x_3x_4^2)dx_1+x_1x_4^2dx_3+(x_4^2-x_1^3)dx_4,$
\item ${\epsilon}_{48}=dx_4,$
\item ${\epsilon}_{58}=-(3x_1^2+x_3x_4)dx_1-x_1x_4dx_3-(2x_4+x_1x_3)dx_4.$
\end{itemize}

In order to compute $A_{ijk}=[v_i,{\epsilon}_{jk}]$ one has to know $[a\frac{\partial}{\partial x_i}, bdx_j]$, where $i,j\in \{ 1,3,4 \}$.\\

The isomorphism ${\wedge}^2T_X\cong {\Omega}^1_X(2)$ identifies $dx_1$ with $\frac{\partial}{\partial x_3}\wedge \frac{\partial}{\partial x_4}$, $dx_3$ with $-\frac{\partial}{\partial x_1}\wedge \frac{\partial}{\partial x_4}$ and $dx_4$ with $\frac{\partial}{\partial x_1}\wedge \frac{\partial}{\partial x_3}$.\\

By definition of the Schouten bracket \cite{Bondal} one has for $j<k$
$$
\left[ a\frac{\partial}{\partial x_i}, b\frac{\partial}{\partial x_j}\wedge \frac{\partial}{\partial x_k}\right]=a\frac{\partial b}{\partial x_i} \left(\frac{\partial}{\partial x_j}\wedge \frac{\partial}{\partial x_k}\right)-b\frac{\partial a}{\partial x_j} \left(\frac{\partial}{\partial x_i}\wedge \frac{\partial}{\partial x_k}\right)+b\frac{\partial a}{\partial x_k} \left(\frac{\partial}{\partial x_i}\wedge \frac{\partial}{\partial x_j}\right),
$$
if $i\neq j, i\neq k$ and
$$
\left[ a\frac{\partial}{\partial x_i}, b\frac{\partial}{\partial x_j}\wedge \frac{\partial}{\partial x_k}\right]=\left(a\frac{\partial b}{\partial x_i}-b\frac{\partial a}{\partial x_i} \right)\left(\frac{\partial}{\partial x_j}\wedge \frac{\partial}{\partial x_k}\right),
$$
if $i= j$ or $i= k$.\\

These formulas allow one to check the expressions for $A_{ijk}$ stated in the Lemma. {\it QED}\\

{\bf Lemma 9.} {\it $B_{ijkl}$ is totally antisymmetric with respect to its indices. In particular, $B_{ijkl}$ are determined completely by the following values:
\begin{gather*}
B_{0123}=-z_{01}+z_{23},\;\; B_{0124}=-z_{15},\;\; B_{0125}=-z_{25},\;\; B_{0128}=-2z_{03}-z_{28},\\
B_{0134}=-z_{08},\;\; B_{0135}=2z_{12}+z_{35},\;\; B_{0138}=z_{38},\;\; B_{0145}=z_{45},\;\; B_{0148}=-z_{48},\\
B_{0158}=2z_{01}-2z_{58},\;\; B_{0234}=-z_{12}+z_{35},\;\; B_{0235}=z_{05}-z_{22},\;\; B_{0238}=z_{08}+3z_{11},\\
B_{0245}=z_{55},\;\; B_{0248}=-z_{01}-z_{58},\;\; B_{0258}=-z_{02}+2z_{45},\;\; B_{0345}=-z_{00}-2z_{15},\\
B_{0348}=2z_{18},\;\; B_{0358}=3z_{03}+4z_{28},\;\; B_{0458}=z_{12}+z_{35},\;\; B_{1234}=-z_{03}+z_{28},\\
B_{1235}=-3z_{00}-z_{15},\;\; B_{1238}=-z_{18}+z_{33},\;\; B_{1245}=-2z_{05},\;\; B_{1248}=2z_{08}+z_{11},\\
B_{1258}=-3z_{12}-4z_{35},\;\; B_{1345}=z_{01}+z_{58},\;\; B_{1348}=-z_{88},\;\; B_{1358}=z_{13}+2z_{48},\\
B_{1458}=z_{03}+z_{28},\;\; B_{2345}=-2z_{02}-z_{45},\;\; B_{2348}=-2z_{13}+z_{48},\\
B_{2358}=-5z_{01}-z_{23}+2z_{58},\;\; B_{2458}=2z_{00}-z_{15},\;\; B_{3458}=-z_{08}+2z_{11}.
\end{gather*}
}

{\it Proof:} We use the same set up and notation as in the proof of Lemma 8 and also restrict to the open set $Z_8=1$.\\

By Proposition 1,
\begin{multline*}
B_{ijkl}=\frac{1}{2}[{\epsilon}_{ij},{\epsilon}_{kl}]=\\
=\frac{1}{2}[x_jdx_i-x_idx_j,x_ldx_k-x_kdx_l]=\frac{1/2}{Vol}\left( (x_jdx_i-x_idx_j)\wedge (dx_l\wedge dx_k-dx_k\wedge  dx_l)+\right.\\
+\left.(dx_j\wedge dx_i-dx_i\wedge  dx_j)\wedge (x_ldx_k-x_kdx_l) \right)=\frac{1}{Vol}\left(  x_i\cdot dx_j\wedge dx_k\wedge dx_l-\right.\\
\left. -x_j\cdot dx_i\wedge dx_k\wedge dx_l+x_k\cdot dx_i\wedge dx_j\wedge dx_l-x_l\cdot dx_i\wedge dx_j\wedge dx_k \right).
\end{multline*}

We can take $Vol=dx_1\wedge dx_3\wedge dx_4$. Then one obtains the expressions stated in the Lemma. {\it QED}\\

\subsection{Variety of Poisson structures.}

Let $\omega=\sum_{i<j} a_{ij}{\epsilon}_{ij}$ be a point in ${\mathbb P}(H^0(X,{\wedge}^2T_X))\cong {\mathbb P}(\mathfrak{so}(7))$. Let us find the equations of the variety of Poisson structures $\mathcal P\subset {\mathbb P}(H^0(X,{\wedge}^2T_X))$ on $X$ in terms of the homogeneous coordinates $a_{ij}$.\\

$$
[\omega,\omega]=\sum_{i<j}\sum_{k<l}a_{ij}a_{kl}[{\epsilon}_{ij},{\epsilon}_{kl}]=4\cdot \sum_{i<j<k<l} {\alpha}_{ijkl}\cdot B_{ijkl},
$$
where ${\alpha}_{ijkl}=a_{ij}a_{kl}-a_{ik}a_{jl}+a_{il}a_{jk}$.\\

Using Lemma 9 one can find the components of $[\omega,\omega]\in H^0(X,{\mathcal O}_X(2))$ relative to the basis $z_{ij}$, $i\leq j$ from Lemma 7. Then $[\omega,\omega]=0$ is equivalent to the following equations:
\begin{gather*}
{\alpha}_{0348}={\alpha}_{0125}={\alpha}_{0138}={\alpha}_{1348}={\alpha}_{0245}={\alpha}_{1238}={\alpha}_{0235}={\alpha}_{1245}=0,\\
{\alpha}_{0123}={\alpha}_{2358},\;\; {\alpha}_{0145}=5{\alpha}_{2345},\;\; {\alpha}_{0148}=5{\alpha}_{2348},\;\; {\alpha}_{0258}=-2 {\alpha}_{2345},\\
{\alpha}_{1358}=2{\alpha}_{2348},\;\; {\alpha}_{0345}+3{\alpha}_{1235}=2{\alpha}_{2458},\;\; {\alpha}_{0248}+6{\alpha}_{2358}=2{\alpha}_{0158}+{\alpha}_{1345},\\
2{\alpha}_{0128}+{\alpha}_{1234}=3{\alpha}_{0358}+{\alpha}_{1458},\;\; {\alpha}_{0134}-{\alpha}_{0238}-2{\alpha}_{1248}+{\alpha}_{3458}=0,\\
-2{\alpha}_{0135}+{\alpha}_{0234}-{\alpha}_{0458}+3{\alpha}_{1258}=0,\;\; {\alpha}_{0124}+2{\alpha}_{0345}+{\alpha}_{1235}+{\alpha}_{2458}=0,\\
{\alpha}_{0128}=4{\alpha}_{0358}+{\alpha}_{1234}+{\alpha}_{1458},\;\; {\alpha}_{0135}+{\alpha}_{0234}+{\alpha}_{0458}-4{\alpha}_{1258}=0,\\
{\alpha}_{0248}+2{\alpha}_{0158}-{\alpha}_{1345}-2{\alpha}_{2358}=0,\;\; {\alpha}_{1248}=-3{\alpha}_{0238}-2 {\alpha}_{3458}.
\end{gather*}

Note that the equations ${\alpha}_{ijkl}=0$, $i<j<k<l$ define the Grassmannian $G(2,7)\subset {\mathbb P}(\mathfrak{so}(7))$ embedded by the Pl{\" u}cker embedding. Provided that one knows that the variety of Poisson structures $\mathcal P$ has two irreducible components of dimensions $10$ and $1$ (this fact is stated in \cite{Loray1}), one concludes immediately that the $10$-dimensional component is exactly this $G(2,7)$, because $dim(G(2,7))=10$. This is, of course, proved more generally in \cite{Loray2}.\\

{\bf Theorem 3 (Loray, Pereira, Touzet, \cite{Loray1}, \cite{Loray2}).} {\it Let $X$ be the (smooth) del Pezzo quintic threefold. Then the variety of Poisson structures $\mathcal P\subset {\mathbb P}(H^0(X,{\wedge}^2T_X))$ on $X$ is the disjoint union of the Grassmannian $G(2,7)\subset {\mathbb P}(H^0(X,{\wedge}^2T_X))$ (embedded by the Pl{\" u}cker embedding) and a smooth conic in ${\mathbb P}(H^0(X,{\wedge}^2T_X))$. The plane spanned by the conic does not intersect the Grassmannian.}\\

{\it Proof:} Let us take the plane $\Pi\subset {\mathbb P}(H^0(X,{\wedge}^2T_X))\cong {\mathbb P}(\mathfrak{so}(7))$ defined by the following linear equations:
\begin{gather*}
a_{02}=a_{05}=a_{08}=a_{13}=a_{15}=a_{18}=a_{24}=a_{25}=a_{34}=a_{38}=a_{45}=a_{48}=0,\\
a_{01}=\frac{5}{2}a_{23},\;\; a_{58}=\frac{9}{2}a_{23},\;\; a_{12}=\frac{5}{3}a_{35},\;\; a_{04}=5a_{35},\;\; a_{03}=\frac{5}{3}a_{28},\;\; a_{14}=-5a_{28}.
\end{gather*}

Then the intersection $\Pi \cap \mathcal P$ is the conic given by the equation $(a_{23})^2=\frac{8}{9}a_{28}a_{35}$ in this plane.\\

Since ${\alpha}_{2358}=a_{23}a_{58}+a_{28}a_{35}=\frac{45}{8}(a_{23})^2$, ${\alpha}_{0345}=-a_{04}a_{35}=-5(a_{35})^2$ and ${\alpha}_{0134}=-a_{03}a_{14}=\frac{25}{3}(a_{28})^2$ never simultaneously vanish on the conic, we conclude that the plane $\Pi$ spanned by the conic does not intersect the Grassmannian $G(2,7)\subset {\mathbb P}(\mathfrak{so}(7))$. {\it QED}\\

\subsection{Poisson cohomology.}

Since $H^i(X,{\wedge}^jT_X)=0$ for any $i\geq 1$ for any $j$ by Lemma 4, it follows from \cite{Stienon} (see Lemma 3.3 in \cite{HongXu}) that Poisson cohomology of $X$ with respect to $\omega\in \mathcal P\subset {\mathbb P}(H^0(X,{\wedge}^2T_X))$ is the cohomology of the following complex:
$$
0 \rightarrow H^0(X,{\mathcal O}_X) \xrightarrow{d_{\omega}=0} H^0(X,T_X) \xrightarrow{d_{\omega}={\alpha}_{\omega}} H^0(X,{\wedge}^2T_X) \xrightarrow{d_{\omega}={\beta}_{\omega}} H^0(X,{\mathcal O}_X(2)) \rightarrow 0.
$$

This implies the following Theorem.\\

{\bf Theorem 4.} {\it Let $X$ be the (smooth) del Pezzo quintic threefold and $\omega\in \mathcal P$ a point on the variety of Poisson structures on $X$. Then
\begin{itemize}
\item $H^0_{Poisson}(X,\omega)\cong H^0(X,{\mathcal O}_X)=k$,
\item $H^1_{Poisson}(X,\omega)\cong ker({\alpha}_{\omega})$,
\item $H^2_{Poisson}(X,\omega)\cong ker({\beta}_{\omega})/ im({\alpha}_{\omega})$,
\item $H^3_{Poisson}(X,\omega)\cong H^0(X,{\mathcal O}_X(2))/ im({\beta}_{\omega})$.
\end{itemize}}

In particular, using matrices $A_{\omega}$ and $B_{\omega}$ computed in Lemma 8 and Lemma 9 we have:
\begin{itemize}
\item $dim(H^0_{Poisson}(X,\omega))=1$,
\item $dim(H^1_{Poisson}(X,\omega))=3-rk(A_{\omega})$,
\item $dim(H^2_{Poisson}(X,\omega))=dim(ker({\beta}_{\omega}))-rk({\alpha}_{\omega})=21-rk(A_{\omega})-rk(B_{\omega})$,
\item $dim(H^3_{Poisson}(X,\omega))=23-rk(B_{\omega})$.
\end{itemize}

\section{Acknowledgement.}
This project started after we read the paper \cite{HongXu}. We thank Jorge Pereira for comments.

\bibliographystyle{ams-plain}

\bibliography{Poisson}

\end{document}